\documentclass{amsart}
\usepackage{amssymb,amsmath}
\usepackage{amsthm}
\usepackage{euscript}
\input{epsf}

\newcounter{sec}

\newcounter{punct}[sec]

\def\punct{\refstepcounter{punct}{\arabic{section}.\arabic{punct}.  }}

\newtheorem{theorem}{Theorem}[sec]
\newtheorem{proposition}[theorem]{Proposition}

\newtheorem{lemma}[theorem]{Lemma}

\newtheorem{corollary}[theorem]{Corollary}

\def\COUNTERS{\addtocounter{sec}{1}
              \setcounter{punct}{0}
          \setcounter{equation}{0}
          \setcounter{theorem}{0}
         }
          
          \def\sm{\smallskip}

\newcommand{\rk}{\mathop {\mathrm {rk}}\nolimits}

\newcommand{\tr}{\mathop {\mathrm {tr}}\nolimits}
\newcommand{\Tr}{\mathop {\mathrm {Tr}}\nolimits}

\def\0{\mathbf 0}

\def\ov{\overline}
\def\wh{\widehat}
\def\wt{\widetilde}

\renewcommand{\rk}{\mathop {\mathrm {rk}}\nolimits}

\newcommand{\PB}{\mathop {\mathrm {PB}}\nolimits}

\def\frg{\mathfrak g}

\def\frc{\mathfrak c}
\def\frd{\mathfrak d}

\def\R {{\mathbb R }}
 \def\C {{\mathbb C }}
  \def\Z{{\mathbb Z}}
  
\def\K{{\mathbb K}}

\def\Q{{\mathbb Q}}

\def\O {{\mathbb O }}

\def\cL{\mathscr L}

\def\cA{\EuScript A}
\def\cO{\EuScript O}

\def\bbO{\mathbb O}

\def\kappa{\varkappa}
\def\epsilon{\varepsilon}
\def\phi{\varphi}
\def\le{\leqslant}
\def\ge{\geqslant}

\def\la{\langle}
\def\ra{\rangle}

\def\O{\mathrm{O}}

\def\U{\mathrm{U}}
\def\GL{\operatorname{GL}}
\def\Sp{\operatorname{Sp}}
\def\End{\operatorname{End}}

\def\dom{\operatorname{dom}}

\def\im{\operatorname{im}}

\def\Assoc{\operatorname{Assoc}}

\def\F{\mathbb{F}}

\def\cL{\EuScript{L}}

\def\cE{\EuScript{E}}

\def\0{\boldsymbol{0}}

\begin{document}

\begin{center}
\Large\bf
Algebra of double cosets of a symmetric group by a smaller symmetric group

\bigskip

\large\sc Yu. A. Neretin%
\footnote{Supported by the grant FWF (The Austrian Scientific Funds),
Project PAT5335224.}

\end{center}

{\small Fix a natural $\alpha$. Let $n\ge \alpha$ be an integer. Consider the symmetric group $S_{\alpha+n}$ and its subgroup $S_n$. We consider the group algebra of  $S_{\alpha+n}$ and its subalgebra $\mathbb{O}[\alpha;n]$ consisting of $S_n$-biinvariant functions, i.e., functions, which are constant on double cosets of $S_{\alpha+n}$ with respect to $S_n$. We obtain two simple descriptions of $\mathbb{O}[\alpha;n]$. First, we write explicitly formulas for multiplication in a natural basis (structure constants  are  Pochhammer symbols). Secondly, we describe this algebra   in terms  of generators and relations. We also construct an interpolating family of algebras $\mathbb{O}[\alpha;\nu]$ depending on a complex parameter $\nu$. 
}

\section{\bf Introduction. Algebras of double cosets}

\COUNTERS

{\bf \punct 
Algebras of double cosets. Preliminaries.%
\label{ss:preliminaries}} For a group $G$ and its subgroup $K$
we denote by $K g K$, where $g\in G$, double cosets of $G$ with respect
to $K$. By $K\backslash G/K$ we denote the space of double cosets,
i.e., the quotient of $G$ with respect to the equivalence relation
$$
 k_1 gk_2\sim g,\qquad\text{where $k_1$, $k_2\in K$.} 
$$

Let a group $G$ be finite. Let $\C[G]$ be its group algebra, we denote
the product in $\C[G]$ (the convolution) as the usual product
by `dot'
(i.e., $\Phi\cdot \Psi$).
 We say, that a function on $G$ is {\it $K$-biinvariant}
 if it is invariant with respect to left and right shifts 
 by elements of $K$.
Equivalently, this function  is constant on double cosets $K g K$.
 Denote 
by
$$
\C[K\backslash G/K]\subset \C[G]
$$
the subalgebra consisting of $K$-biinvariant functions
(for more details, see, e.g., \cite{Bum}).

For $g\in G$, denote by $\delta_g$ the element $g$ considered as an 
element of $\C[G]$. Denote 
$$
\delta_K=\frac{1}{\#K} \sum_{k\in K}\delta_K,
$$
where the symbol $\# X$ denotes the number of elements of a finite set $X$. 
 For a double coset $\frg=KgK$,  we define
\begin{equation}
 \delta_\frg=\delta_K\cdot\delta_g\cdot\delta_K
=\frac{1}{\# (KgK)}\sum_{q\in KgK} \delta_q.
\label{eq:delta-frg}
 \end{equation}
Clearly, the elements $\delta_\frg$, where $\frg$ ranges in $K\backslash G/K$, form a basis
 in the algebra $\C[K\backslash G/K]$. 
 
 Let $\rho$ be an irreducible representation of $G$
 in a Euclidean space $H_\rho$. Denote by $H_\rho^K\subset H_\rho$
 the subspace of $K$-fixed vectors (it can be trivial).
 Denote by $P_\rho^K$ the operator
 of orthogonal projection to $H_\rho^K$,
 $$
 P_\rho^K=\rho(\delta_K).
 $$
 For $g\in G$ we define the operator
 $$
 \wt\rho(g):H_\rho^K \to H_\rho^K
$$
by 
$$
\wt \rho(g):=P_\rho^K\rho(g)\Bigr|_{H_\rho^K}.
$$
It is easy to see that this operator depends only on the double
coset $\frg$ containing $g$, so we can write $\wt\rho^K(\frg)$. Decomposing%
\footnote{The symbol $\bot$ denotes an orthogonal complement.}
$H_\rho=H^K_\rho\oplus(H_\rho^K)^\bot$, we get
$$
\rho(\delta_\frg)=\rho(\delta_K)\,\rho(g)\,\rho(\delta_K)
=\begin{pmatrix}
\wt \rho^K(\frg)&0\\0&0
\end{pmatrix}.
$$ 
So, $\wt\rho(\frg)$ is a representation of the algebra
$\C[K\backslash G/K]$ in $H_\rho^K$.

This easily implies that $\C[K\backslash G/K]$ is isomorphic
to a direct sum
of matrix algebras,
\begin{equation}
\C[K\backslash G/K]\simeq \bigoplus_{\rho\in \wh G}
\End(H_\rho^K),
\label{eq:oplus} 
\end{equation}
where the summation is taken over the set $\wh G$ of
 all irreducible representations 
of $G$ (actually we have only summands, for which
$H_\rho^K\ne 0$).
So,
$$
\dim \C[K\backslash G/K]= \sum_{\rho\in \wh G}
(\dim H_\rho^K)^2.
$$

\sm

Next, we define some additional simple structures
on algebras $\C[K\backslash G/K]$.

\sm

1) We have a {\it trace} on $\C[G]$ defined by
\begin{equation*}
\tr \Bigl(\sum_{g\in G} c_g \delta_g\Bigr):= c_1,
\end{equation*}
where $1$ is the unit of the group. It is easy to see
that this trace satisfies 
$$
\tr (\Phi\cdot\Psi)=\tr(\Psi\cdot\Phi).
$$
Restricting the trace to $\C[K\backslash G/ K]$
we get a trace on $\C[K\backslash G/ K]$,
\begin{equation}
\tr \Bigl(\sum_{\frg\in K\backslash G/ K} a_\frg\delta_\frg \Bigr)=
\frac 1{\# K} a_K. 
\label{eq:trace}
\end{equation}

\sm

2) We  define an antilinear {\it involution} in $\C[G]$ by
$$
\bigl(\sum c_g \delta_g\bigr)^*:=\sum \ov c_g \delta_{g^{-1}}.
$$ 
Clearly,
$$
(\Phi+\Psi)^*=\Phi^*+\Psi^*,\quad (\Phi\cdot\Psi)^*=\Psi^*\cdot\Phi^*,
\quad \Phi^{**}=\Phi \quad (s\, \Phi)^*=\ov s \,\Phi^*,
$$
where $\Phi$, $\Psi\in\C[G]$, $s\in \C$.
Clearly, $\C[K\backslash G/K]$ is invariant with respect to the involution; 
$$
\Bigl(\delta_{KgK}\Bigr)^*=\delta_{Kg^{-1}K}.
$$

3) We  have a positive definite  {\it inner product}
$$
\la \Phi,\Psi\ra:=\tr (\Phi\cdot \Psi^*),
$$
Equivalently,
$$
\bigl\la \sum b_g \delta_g, \sum c_g \delta_g\bigr\ra
=\sum_{g\in G} b_g \ov c_g.
$$
So, we get  the inner product on $\C[K\backslash G/K]$,
$$
\bigl\la \sum b_g \delta_\frg, \sum c_\frg \delta_g\bigr\ra
=\sum_{\frg\in K\backslash G/K} \frac{1}{\#\frg}\, b_\frg \ov c_\frg.
$$

4) We have a homomorphism $\iota:\C[G|\to \C$ defined by
\begin{equation}
\iota \Bigl(\sum c_g \delta_g\Bigr):=\sum c_g
\label{eq:iota}
\end{equation}
(we take the trivial representation of $G$ in $\C$, the corresponding
representation of the group algebra is $\iota$). 
The restriction of $\iota$
 to $\C[K\backslash G/K]$ is given by
$$
\iota\Bigl(\sum_{\frg\in K\backslash G/ K} a_\frg\delta_\frg \Bigr)=
\sum a_\frg.
$$

\sm

These considerations
%\footnote{With some care with  traces.} 
can be easily extended to the case, when 
$G$ is a locally compact group and $K$ is a compact subgroup.
Now the convolution algebra $\C[K\backslash G/K]$
consists of (compactly supported) $K$-biinvariant complex-valued measures
on $G$ (here we can take different functional spaces consisting of
$K$-biinvariant functions or distributions, such space must be closed
with respect to the convolution).

\sm

{\bf \punct Some important cases.}

\sm

a) Let $G$ be a linear semisimple Lie group and $K\subset G$ the maximal 
compact subgroup.
According  Gelfand theorem  \cite{Gel},
the algebras $\C[K\backslash G/K]$ are commutative, for this reason
$\dim H_\rho^K\le 1$. The same statement holds if $G$ is a compact Lie group and $K$ is its symmetric subgroup.
 These algebras were
discussed and used a lot starting \cite{Gel} and \cite{BG}; as far as I know their nice transparent descriptions are
not invented. 

\sm

Consider the case of rank one simple Lie groups.
Let $\K$ denote the reals $\R$, complex numbers $\C$, 
or quaternions $\mathbb{H}$, denote $\kappa:=\dim\K$. Consider a space
$\K^n$
 equipped with an Hermitian form with inertia indices
$(n-1,1)$, i. e., the form defined by 
the $((n-1)+1)$-block matrix $J=\begin{pmatrix}1&0\\0&-1\end{pmatrix}$.
 Consider the groups $G:=\U(n-1,1;\K)$
preserving these forms (a matrix $g=\begin{pmatrix}a&b\\c&d\end{pmatrix}$ over $\K$ is contained in $G$
if $gJg^*=J$).
Their maximal compact subgroups are $K=\U(n-1;\K)\times \U(1;\K)$
(they consists of matrices  $g=\begin{pmatrix}a&0\\0&d\end{pmatrix}$
such that $aa^*=1$, $d\,\ov d=1$).

In the standard notation, for $\K=\R$, $\C$, $\mathbb{H}$  a group 
$G=\U(n-1,1;\K)$ is respectively
\begin{equation}
G=\O(n-1,1),\quad \U(n-1,1), \quad \Sp(n-1,1),
\label{eq:ranl-one}
\end{equation}
and $K$ is 
$$
K=\O(n-1)\times \O(1),\quad \U(n-1)\times \U(1), \quad \Sp(n-1)\times \Sp(1).
$$
In all cases,  convolutions of  delta-functions
$\delta_{KgK}$
are determined by  ${\vphantom{F}}_2F_1$-hy\-per\-ge\-omet\-ric expressions.
There is a common formula for these expressions including $n$
and $\kappa$ as parameters. More precisely, 
there is a meromorphic family of associative algebras depending on two complex parameters
$n$, $\kappa$ whose elements are
functions on the half-line, and a product
$$
f_1\,\bullet_{\vphantom{\bigl|}n,\kappa}\, f_2(z)=\int_0^\infty \cL_{n,\kappa}(x,y;z)\, f_1(x)\,f_2(y)\,
dx\,dy
$$
 is determined by a certain hypergeometric kernel $\cL_{n,\kappa}(z;x,y)$,
  which  for  integer positive $n$
and $\kappa=1$, $2$, $4$ gives algebras $\C[K\backslash G/K]$ 
for the cases \eqref{eq:ranl-one} (for $n=2$ and $\kappa=8$
we get the double coset algebra for one of real forms 
of the exceptional simple Lie group $F_4$, the number 8 is the dimension of the octaves). See Flensted-Jensen and Koornwinder \cite{FJK},  \cite{Koo}, Sect 7.
\hfill $\boxtimes$

\sm

2) {\it Hecke algebras}, Iwahori \cite{Iwa}. Consider a group $G:=\GL(n;\F_q)$
 over a finite field $\F_q$ with $q$ elements,
let $K$ be the group of upper-triangular matrices. It is easy to see that 
each
double coset $K\backslash G/K$ contains a unique matrix-permutation.
Denote by $T_i$, where $1\le i\le n-1$, the double coset containing the matrix of 
the transposition $\bigl(i(i+1)\bigr)$. Then $\C[K\backslash G/K]$
is the algebra with generators $T_i$ and relations
$$
(T_i-q)(T_i+1)=0,\quad
T_i T_{i+1}T_i=T_{i+1} T_i T_{i+1},\quad
T_i T_j=T_j T_i \quad\text{if $|i-j|\ge 2$.}
$$ 
This algebra makes sense for all $q\in \C$ 
and always has dimension $n!$. For $q=1$
we get the group algebra of the symmetric group $S_n$
(indeed, the first relation gives $T_i^2=1$, therefore $(T_i T_{i+1})^3=1$,
and we get relations for $S_n$ generated by transpositions). See,
e.g., \cite{DG}, \cite{Bum}, Sect.~48. This algebra can live by its own life
without the groups $\GL(n;\F_q)$ and $S_n$.

If $G:=\GL(n;\F_q)$ and $K$ is the group of strictly upper triangular
matrices,
then the algebra of double cosets admits a nice description
(Yokonuma \cite{Yok}), but it does not admit an interpolation in the parameter $q$.
\hfill $\boxtimes$

\sm

3) {\it The  affine Hecke algebra}, Iwahori, Matsumoto \cite{IwaMats}.
Let $\Q_p$ be $p$-adic field, $\Z_p$ be $p$-adic integers.
 Now $G$ is the general  group $\GL(n;\Q_p)$
over a $p$-adic field, $K$ is the Iwahori subgroup, i.e., the group of matrices
whose entries

---  over the diagonal $\in \Z_p$, 

---  on the diagonal
$\in\Z_p\setminus p\Z_p$,

---  and under the diagonal $\in p\Z_p$.

 Then $\C[K\backslash G/K]$ has a nice description and admits an interpolation with respect to $p$.
Again, this algebra lives by its own life independently on groups
over locally compact fields, see \cite{Mac1}. On the variant
$G=\GL(n,\Q_p)$, $K=\GL(n,\Z_p)$, see Maconald \cite{Mac}, Chapter V.
\hfill $\boxtimes$

\sm

Algebras of types 1)-3) have some relatives, but in any case
a  list of known  algebras
of double cosets with simple nice descriptions is short%
\footnote{Let $G$ be a symmetric group $S_{m+n}$ permuting 
the set $\{1,\dots,n+m\}$,
$K$ be $S_n\times S_m$ embedded to $S_{n+m}$ in the natural way.
Double cosets  are enumerated by a number of $i\le n$ such that $g(i)>n$.
Structure constants in the natural basis
\eqref{eq:delta-frg} are hypergeometric
expressions of type $\vphantom{F}_4F_3$ at 1 (see \cite{Ner-zametki}).
Such  description does not look as too opthimistic
 since it is  one simplest
double coset spaces $K\backslash G/K$. In any case, this example indicates
difficulties related to our topic. On the other hand, the condition of associativity of $\C[K\backslash G/K]$
implies unusual quadratic identities for $\vphantom{F}_4F_3[1]$.}.
More general algebras of double cosets were topics
of numerous efforts, see, e.g., \cite{FH}, \cite{IK}, \cite{Mel}, \cite{Tou}, \cite{KR}.

\sm

For infinite-dimensional groups $G$
quite often there appear natural associative multiplications on 
sets $K\backslash G/K$, such semigroups act in spaces of $K$-fixed
vectors of unitary representations of $G$. These semigroups admit
transparent descriptions, see, e.g.,
 \cite{Olsh1}, \cite{Olsh}, \cite{Ner-bist},
  \cite{Ner-book}, Chapter VIII, Sect.IX.3-4, \cite{Ner-umn}, \cite{Ner-coll}, \cite{Ner-finite}.
 There arises a question about descent to finite-dimensional level.   The author did several experiments in this direction, see, e.g.,
 \cite{Ner-zametki}, \cite{Ner-ivanov} (and  \cite{Ner-PBL},
 which is a partial extension of the present work to groups 
 $\GL(n,\F_q)$ over finite fields). 
 
\sm

 {\bf \punct The purpose of the paper.} Fix $\alpha>0$, $n\ge \alpha$.
 Consider the symmetric group $G=S_{\alpha+n}$.
 Double cosets $S_n\backslash S_{\alpha+n}/S_n$
 are enumerated by partial  bijections of the set $\{1,2,\dots,\alpha\}$.
 %For algebras of 
 We consider algebras of double cosets 
 $$\bbO[\alpha;n]:=\C[S_n\backslash S_{\alpha+n}/S_n] 
 \quad\text{for $n\ge\alpha$.}$$
 For these algebras:

\sm

--- we evaluate the product in the natural basis $\delta_\frg$,
 the formula for structure constants is unexpectedly simple (Sect.~2);

\sm 
 
---  we  construct a family $\bbO[\alpha;\nu]$
 of algebras depending on a complex parameter $\nu$
 interpolating algebras $\C[S_n\backslash S_{\alpha+n}/S_n]$ (Sect.~3);

\sm

---  we obtain a   description $\bbO[\alpha;\nu]$ in  terms of generators
 and relations (Sect.~4).

\section{Structure constants}

\COUNTERS

This section contains preliminaries on partial bijections
(Subsect. \ref{ss:partial-bijections}--\ref{ss:corners}),
 description of double cosets
$S_n\backslash S_{\alpha+n}/S_n$ (Subsect. \ref{ss:cosets}), 
evaluation of structure constants
of $\bbO[\alpha;n]$ (Theorem \ref{th:coefficients}),
decomposition of $\bbO[\alpha;n]$ to a direct sum (Proposition
\ref{pr:structure}).

\sm

{\bf \punct Some notation.}
For integer $\alpha\ge 0$ denote by $\Omega_\alpha$
 the set $\{1,2,\dots,\alpha\}$.

We say that a {\it $0$-$1$-matrix} is a matrix composed of zeros and units and
containing $\le 1$ unit in each column and in each row.
A square 0-1-matrix is invertible if and only if it contains precisely
 one unit in each column and each row.

 We denote by $S_m$ the {\it symmetric group}
of order $m$, i.e., the group of all permutations of the
set $\Omega_m$.   We also regard a symmetric group as the group of 
of all invertible 0-1-matrices of order $m$.

 Denote by $(ij)\in S_m$ the transposition of $i$ and $j$.

\sm

{\bf \punct Partial bijections.%
\label{ss:partial-bijections}}
Let $A$, $B$ be a finite sets. We say that a {\it partial bijection}
$\sigma:A\to B$ is a bijective map from a subset 
$X\subset A$ to a subset $Y\subset B$, we
say that $X$ is the {\it domain} $\dom\sigma$ of $\sigma$, $Y$ is the {\it image} $\im\sigma$
 of $\sigma$.
The {\it rank} of a partial bijection $\sigma$ is
 $$
 \rk \sigma:=\# \dom \sigma= \# \im \sigma
 .$$ 
 We define the {\it pseudoinverse} partial bijection
$\sigma^\square:B\to A$ with $\dom \sigma^\square=\im\sigma$,
$\im \sigma^\square=\dom\sigma$, and $\sigma^\square(y)=x$
iff $\sigma(x)=y$.
 
If $\tau:A\to B$, $\sigma: B\to C$ are partial bijections, 
then their {\it product} $\sigma\tau:A\to C$ is
a natural product of partially defined maps.
Precisely, $\sigma\tau$
  is defined
from  the following conditions:

\sm

---
$
\dom \sigma\tau$ is the set of all $a\in \dom \tau$ such that $\tau(a)\in \dom \sigma$;

\sm

---  for any $a\in \dom \sigma\tau$ we have 
 $(\sigma\tau)(a)=\sigma(\tau(a))$.

\sm

Clearly, $(\sigma \tau)^\square=\tau^\square\sigma^\square$.

\sm

Obviously,
$$
\dom\sigma\tau\subset \dom\tau,\quad
\im\sigma\tau\subset \im\sigma,\quad \rk \sigma\tau\le\min(\rk\sigma,\rk\tau).
$$

We denote  by
 $\PB(\alpha)$ the semigroup of all partial bijections of 
$\Omega_\alpha:=\{1,\dots,\alpha\}$ equipped with the involution
$\sigma\to \sigma^\square$. The identical map is the unit of
$\PB(\alpha)$, the group of invertible elements is
 $S_\alpha$.
 
 \sm

For a partial bijection 
$\sigma:\Omega_\alpha\to\Omega_\beta$
we define the $0$-$1$-matrix  $S=\{s_{ij}\}$ of size $\beta\times\alpha$
assuming $s_{ij}=1$ if and only if $\sigma(j)=i$.

Clearly, a product of partial bijections corresponds to a product 
of $0$-$1$-matrices, $\rk S=\rk\sigma$, and the pesudoinversion
corresponds to the transposition. {\it Below we do not distinguish
a partial bijection and the corresponding $0$-$1$-matrix.}

\sm

For a subset $\Delta\subset \Omega_\alpha$
denote by $T(\Delta)\in\PB(\alpha)$ the partial
bijection such that  

\sm

--- $\dom\bigl(T(\Delta)\bigr)=\Delta$;

\sm

--- $T(\Delta)$ is identical on $\Delta$. 

\sm

Clearly, these elements are {\it idempotents}, $T(\Delta)^2=T(\Delta)$.
Moreover,
$$
T(\Delta_1)\,T(\Delta_2)= T(\Delta_1\cap \Delta_2)=T(\Delta_2)\,T(\Delta_1).
$$

Obviously, any $\sigma\in \PB(\alpha)$ can be represented in forms 
\begin{align}
\sigma&= g\, T(\dom\sigma), \qquad\text{where $g\in S_\alpha$;}
\label{eq:gT}
\\
\sigma&=T(\im\sigma)\, h, \qquad\text{where $h\in S_\alpha$.}
\label{eq:Tg}
\end{align}

\sm

Let $\phi$, $\psi\in \PB(\alpha)$.
We say that $\psi$ is an {\it extension} of $\phi$,
if $\dom\psi\supset \dom\phi$ and the restriction
 of $\psi$ to $\dom\phi$ coincides with $\phi$ 
 (in particular, $\phi$ is an extension of itself). Equivalently,
 a $0$-$1$-matrix of $\phi$ can be obtained by 
 replacing of some units by zero. 
We denote 
$$
\psi\sqsupset \phi,\quad\text{or}\quad \phi\sqsubset \psi. 
$$

\begin{lemma}
\label{l:number-bijections}
The number of partial  bijections 
$\Omega_\alpha\to \Omega_\beta$ of rank $\rho$ is
\begin{equation}
\frac{\alpha!\,\beta!}{\rho!\,(\alpha-\rho)!\,(\beta-\rho)!}
.
\label{eq:number-rank}
\end{equation}
\end{lemma}

{\sc Proof.}
The group $S_\beta\times S_\alpha$ acts on the set $U_\rho$
of such partial bijections  by left
and right multiplications. 
Clearly, $U_\rho$ is homogeneous.
The stabilizer of the idempotent
$T(\Omega_\rho)$
 consists of pairs of matrices
of the form
$$
\begin{pmatrix}
h&0\\0&p
\end{pmatrix}\in S_\beta=S_{\rho+(\beta-\rho)},
\qquad
\begin{pmatrix}
h^{-1}&0\\0&q
\end{pmatrix}\in S_\alpha=S_{\rho+(\alpha-\rho)}.
$$
Hence, the stabilizer
is isomorphic to the group 
$S_\rho\times S_{\beta-\rho}\times S_{\alpha-\rho}$.
Thus,
$$
\#U_\rho=\frac{\#(S_\alpha\times S_\beta}{\# (S_\rho\times S_{\beta-\rho}\times S_{\alpha-\rho})},
$$
i.e., \eqref{eq:number-rank}. 
\hfill $\square$

\sm 

{\sc Remark.} Therefore, we get
\begin{multline*}
\# \PB(\alpha)
=
\sum_{\rho=0}^\alpha
\frac{(\alpha!)^2}{((\alpha-\rho)!)^2\,\rho!}
=\sum_{r=0}^\alpha
\frac{(\alpha!)^2}{(\alpha-r)!\,(r!)^2}
=\\=
\alpha!\sum_{r=0}^\alpha \frac{(-\alpha)_r}{(r!)^2}(-1)^r
=\alpha!\, \vphantom{F}_1F_1[-\alpha;1;-1],
\end{multline*}
where
$$
(a)_k:=a(a+1)\dots(a+k-1)
$$
is the {\it Pochhammer symbol} and
 $\vphantom{F}_1F_1[\cdot]$ denotes the Kummer confluent hypergeometric function. \hfill $\boxtimes$
 
 \sm

{\bf\punct Corners of matrices-permutations.%
\label{ss:corners}}
For a matrix $Q$ of size $M\times M$ denote by $[Q]_{\mu\nu}$
its left upper corner of size $\mu\times \nu$. 
%For $Q\in S_{n+\alpha}$ we denote by $[Q]_\alpha$ 
%its left upper corner of size $\alpha$.
So the map $Q\mapsto [Q]_{\mu\nu}$
sends elements of the symmetric group $S_M$
to partial bijections $\Omega_\nu\to \Omega_\mu$.

Clearly, 
$$
M-\mu-\nu+\rk\, [Q]_{\mu\nu}\ge 0.
$$

\begin{lemma}
\label{l:number-preimage}
Let $\sigma$ be a partial bijection  $\Omega_\alpha\to \Omega_\beta$,
 $\rk(\sigma)=\rho$. 
Then the number of $Q\in S_M$ such that $[Q]_{\beta\alpha}=\rho$
is
\begin{equation}
\frac{(M-\alpha)!\,(M-\beta)!}{(M-\alpha-\beta+\rho)!}.
\label{eq:number-preimage}
\end{equation}
\end{lemma}

{\sc Proof.} Without lost of generality, we can assume
that $[Q]_{\beta\alpha}$ has the following $(\rho+(\beta-\rho))\times (\rho+(\alpha-\rho))$-block structure:
$$
[Q]_{\beta\alpha}=\begin{pmatrix}
1&0\\0&0
\end{pmatrix}.
$$
Therefore, $Q$ has the form
$$
\left(
\begin{array}{cc|c}
1&0&0\\0&0&S\\
\hline
0&T&R
\end{array}
\right).
$$
The size of a $0$-$1$-matrix $S$ is $(\beta-\rho)\times (M-\alpha)$,
it contains a unit in each  row. The number of such matrices
is
\begin{equation}
(M-\alpha)(M-\alpha-1)\dots (M-\alpha-\beta+\rho+1)
=
\frac{(M-\alpha)!}{(M-\alpha-\beta+\rho)!}.
\label{eq:factor1}
\end{equation}
In the same way, a matrix $T$ can be chosen 
in
\begin{equation}
\frac{(M-\alpha)!}{(M-\alpha-\beta+\rho)!}.
\label{eq:factor2}
\end{equation}
ways. After fixing   $S$ and $T$, we must chose the remaining
$(M-\alpha-\beta+\rho)$
units in $R$. They must be located in $(M-\alpha-\beta+\rho)$
columns and $(M-\alpha-\beta+\rho)$ rows. So, we can chose them 
in $(M-\alpha-\beta+\rho)!$ ways.
Multiplying this by \eqref{eq:factor1} and \eqref{eq:factor2},
we get \eqref{eq:number-preimage}.
\hfill $\square$

\sm

%Lemmas \ref{l:number-bijections} and \ref{l:number-preimage}
%imply the following corollary:

%\begin{corollary}
%The number of elements $Q\in S_M$ such that
%$\rk[Q]_{\beta\alpha=\rho}$
%is 
%$$
%\frac{\alpha!\, (M-\alpha)!\,\beta!\,(M-\beta)!}
%{\rho!\,(\alpha-\rho)!\,(\beta-\rho)!\,(M-\alpha-\beta+\rho)!}
%$$
%\end{corollary}

%{\bf \punct Extensions of partial bijections.}
%Let  $\sigma\in \PB(\alpha,\beta)$. We say 
%that $\kappa \in \PB(\alpha,\beta)$ is an {\it extension} of $\sigma$
%if 
%$$
%\dom\kappa\supset \dom \dom \sigma
%\quad\text{and}\quad \kappa\Bigr|_{\dom \sigma}=\sigma
%$$
%We denote by $\Ext(\sigma)$ the set of all  extensions of $\sigma$

\sm

{\bf\punct  Double cosets 
$\boldsymbol{S_n\backslash S_{\alpha+n}/S_n }$.%
\label{ss:cosets}}
Let $S_{\alpha+n}$ denote the group
of permutations of the set 
$\Omega_{\alpha+n}$.
{\it Below $\alpha$ is fixed and $n$ varies.}
 Let $S_\alpha$ be the group of permutations of 
 the set $\Omega_\alpha=\{1,\dots,\alpha\}$, and $S_n$
 be the group of permutations  of
% the set 
 $\Omega_{\alpha+n}\setminus\Omega_\alpha=\{\alpha+1, \dots, \alpha+n\}$.

%Let $i_1<i_2<\dots < i_k\le \alpha$.
%We denote by
%\begin{equation}
%T_{i_1,\dots,i_k}
%\label{eq:idempotent}
%\end{equation}
%the diagonal element of $\Pi_\alpha$ whose $i_1$, \dots, $i_k$-th
% entries are 0 and the remaining diagonal entries are units.
% Such  elements are idempotents in $\Pi_\alpha$.

For the following  simple statement  see, e.g., \cite{Olsh} 
or \cite{Ner-book}, Sect. VIII.1.

\begin{theorem}
{\rm a)} The map from $S_{\alpha+n}\to \PB(\alpha)$ defined by
$Q\mapsto [Q]_{\alpha\alpha}$
is constant on double cosets $S_n\cdot g\cdot S_n$
and determines an injective map 
\begin{equation}
S_n\backslash S_{\alpha+n}/S_n\to \PB(\alpha).
\label{eq:SSSPB}
\end{equation}

\sm

{\rm b)} If $n\ge \alpha$, then the map \eqref{eq:SSSPB}
 is a bijection.
If $n<\alpha$, then the image consists of elements
of $\PB(\alpha)$ of rank $\ge\alpha-n$.
\end{theorem}

{\bf \punct Structure constants.}
Our topic is  algebras 
$$\bbO[\alpha;n]:=\C[S_n\setminus S_{\alpha+n}/S_n]
\quad
\text{for $n\ge \alpha$.} $$

We consider the following basis in $\C[S_n\backslash S_{\alpha+n}/S_n]$
enumerated by elements%
\footnote{Elements of this basis differs from elements 
\eqref{eq:delta-frg} by scalar factors depending on elements.}
  $\sigma\in \PB(\alpha)$:
$$
\Xi[\sigma]=\frac{1}{n!}
\sum_{g\in S_{\alpha+n}:\,[g]_{\alpha\alpha}=\sigma}
\delta_g.
$$

\begin{theorem}
\label{th:coefficients}
Let $\sigma$, $\tau\in \PB(\alpha)$.
Then
\begin{equation}
\Xi(\sigma)\,\Xi(\tau)=
\sum_{\begin{array}{l}\kappa\sqsupset \sigma\tau:
\\ \dom\kappa\cap
 \dom\tau=\dom\sigma\tau,
 \\
\im\kappa\cap \im \sigma=\im\sigma\tau\end{array}} 
\frc^\kappa_{\sigma, \tau}(n)\, \Xi(\kappa),
\label{eq:coefficients}
\end{equation}
where the structure constants $\frc^\kappa_{\sigma, \tau}(n)$
are given by
\begin{multline}
\frc^\kappa_{\sigma, \tau}(n)=
\frac{(n-\alpha+\rk \kappa)!}
{(n-2\alpha+\rk \sigma+\rk \tau-\rk \sigma\tau+\rk\kappa )!}=
\\=
(n-\alpha+\rk \kappa)(n-\alpha+\rk \kappa-1)\dots 
(n-2\alpha+\rk \sigma+\rk \tau-\rk \sigma\tau+\rk\kappa +1).
\label{eq:constants}
\end{multline}
\end{theorem}

{\sc Remark.} We have
\begin{equation}
\bigl(n-\alpha+\rk \kappa\bigr)-(n-2\alpha+\rk \sigma+\rk \tau-\rk \sigma\tau+\rk\kappa\bigr)=\alpha+\rk \sigma\tau-\rk\sigma -\rk\tau\ge 0,
\label{eq:hvatit}
\end{equation}
i.e., the number of factors in
\eqref{eq:constants} is non-negative.
\hfill $\boxtimes$

\sm

{\bf \punct Proof of Theorem \ref{th:coefficients}.}

\begin{lemma}
Let  $\sigma$, $\tau$ range in $\PB(\alpha)$.
Consider the natural action of the group $S_\alpha\times S_\alpha\times S_\alpha$ 
 on pairs $(\sigma, \tau)\in \PB(\alpha)\times \PB(\alpha)$:
\begin{equation}
(g_1,g_2,g_3):\,(\sigma,\tau)\mapsto
 \bigl(g_1 \sigma g_2^{-1}, g_2 \tau g_3^{-1}\bigr).
\label{eq:ggg}
\end{equation}
Its orbits are enumerated by triples
\begin{equation}
(\rk \sigma,\,\rk \tau,\,\rk \sigma\tau).
\label{eq:triples}
\end{equation}
\end{lemma}

{\sc Proof.} Fix a collection \eqref{eq:triples}.
Clearly, the second copy of $S_n$ acts transitively
on the set of all possible triples 
$(\im \tau, \dom\sigma, \im\tau\cap\dom\sigma)$.
After fixing these sets, the first copy of $S_n$ can move $\sigma$
to an arbitrary desired position (see \eqref{eq:gT}). Also, the third copy can move $\tau$
 to an arbitrary desired position (see \eqref{eq:Tg}).
\hfill  $\square$   

\sm

Define $i$, $j$, $k$ by
$$i:=\rk \sigma \tau,\quad i+j:=\rk\sigma,\quad i+k=\rk \tau,$$
and let $l$ be defined from the condition
\begin{equation}
\alpha=i+j+k+l.
\label{eq:ijkl}
\end{equation}
Notice that $i$, $j$, $k$, $l$ are nonnegative.
 By the lemma, transformations \eqref{eq:ggg}
allow us to reduce 
a pair $(\sigma,\tau)$ to the form
$$
\sigma=T\{1,2,\dots,i,i+1\dots,i+j\},\qquad
\tau=T\{1,2,\dots,i,i+j+1,\dots,i+j+k\},
$$ 
where $T\{\cdot\}$ are idempotents defined in Subsect. \ref{ss:partial-bijections}.
Therefore, it is sufficient to  verify our statement for such cases.

\sm

We take the following representatives of the corresponding double cosets:
$$J_\sigma:=\left(
\begin{array}{cccc|ccccc}
1& 0& 0& 0& 0& 0& 0& 0& 0\\ 
0& 1& 0& 0& 0& 0& 0& 0& 0\\ 
0& 0& 0& 0& 0& 0& 1& 0& 0\\
0& 0& 0& 0& 0& 0& 0& 1& 0\\
\hline
0& 0& 0& 0& 1& 0& 0& 0& 0\\
0& 0& 0& 0& 0& 1& 0& 0& 0\\ 
0& 0& 1& 0& 0& 0& 0& 0& 0\\
0& 0& 0& 1& 0& 0& 0& 0& 0\\
0& 0& 0& 0& 0& 0& 0& 0& 1
\end{array}\right)
\begin{array}{l}
\dots i\\
\dots j\\
\dots k\\
\dots l\\
\dots i\\
\dots j\\
\dots k\\
\dots l\\
\dots n-\alpha
\end{array}
$$
and
$$
J_\tau:=\left(
\begin{array}{cccc|ccccc}
1& 0& 0& 0& 0& 0& 0& 0& 0\\
0& 0& 0& 0& 0& 1& 0& 0& 0\\
0& 0& 1&   0& 0& 0& 0& 0& 0\\
0& 0& 0& 0& 0& 0& 0& 1& 0\\
\hline
0& 0& 0& 0& 1& 0&   0& 0& 0\\
0& 1& 0& 0& 0& 0& 0& 0& 0\\
0& 0& 0& 0& 0& 0& 1& 0&   0\\
0& 0& 0& 1& 0& 0& 0& 0& 0\\
0& 0& 0& 0& 0& 0& 0& 0& 1
\end{array}\right).
$$
We must examine double cosets of the form
$$
S_n\cdot J_\sigma \cdot
\left(\begin{array}{c|c}
1&0\\
\hline
0&D
\end{array}\right)\cdot J_\tau\cdot S_n,
$$
where $D$ ranges in $S_n$. Representing $D$ as a block matrix of size
$(i+j+k+l+(n-\alpha))$ with blocks $\{d_{ij}\}_{1\le i,j\le 5}$ we get
$$
 J_\sigma \cdot
\left(\begin{array}{c|c}
1&0\\
\hline
0&D
\end{array}\right)\cdot J_\tau=
\left(
\begin{array}{cccc|ccccc}
1& 0& 0& 0& 0& 0& 0& 0& 0\\
0& 0& 0& 0& 0& 1& 0& 0& 0\\
0& d_{32}& 0& d_{34}& d_{31}& 0& d_{33}& 0& d_{35}\\
0& d_{42}& 0& d_{44}& d_{41}& 0& d_{43}& 0& d_{45}\\
\hline
0& d_{12}& 0& d_{14}& d_{11}& 0& d_{13}& 0& d_{15}\\
0& d_{22}& 0& d_{24}& d_{21}&   0& d_{23}& 0& d_{25}\\
0& 0& 1& 0& 0& 0& 0& 0& 0\\
0& 0& 0& 0& 0& 0& 0&  1& 0\\
0& d_{52}& 0& d_{54}& d_{51}& 0& d_{53}& 0& d_{55}
\end{array}\right).
$$
So, for a given $D$ we get the following element of $\PB(\alpha)$: 
$$
\kappa_D=\left(\begin{array}{ccccl}1&0&0&0&
\\
0&0&0&0
\\
0&d_{32}&0&d_{34}
\\
0&d_{42}&0&d_{44}
 \end{array}\right) \begin{array}{l}
\dots i\\
\dots j\\
\dots k\\
\dots l
\end{array} .
$$
It satisfies the following properties:

\sm

--- $\kappa$ is an extension of 
$\sigma\tau=\begin{pmatrix}
1&0&0&0\\
0&0&0&0\\
0&0&0&0\\
0&0&0&0
\end{pmatrix}
$;

\sm

--- $\dom \kappa_D\subset \{1,\dots,i\}\cup\{i+1,\dots,i+j\} \cup \{i+j+k+1,\dots,\alpha\}$;

\sm

--- $\im \kappa_D\subset  \{1,\dots,i\}\cup\{i+j+1,\dots,i+j+k\} \cup \{i+j+k+1,\dots,\alpha\}$.

\sm

So, {\it we obtained the set of summation in formula \eqref{eq:coefficients}}.

\sm

Clearly, any $0$-$1$-matrix of size $(k+l)\times(j+l)$
can be represented as a submatrix
$\begin{pmatrix}
d_{32}&d_{34}\\
d_{42}&d_{44}
\end{pmatrix}
$   
of $D=\{d_{ij}\}_{1\le i,j\le 5}$
 of $S_n$. Denote by $\rho$ the rank of this matrix.
 By Lemma \ref{l:number-preimage}, we can   embed it to 
 a matrix of permutation in
\begin{equation}
\frac{(n-k-l)!\,(n-j-l)!}{(n-j-k-2l+\rho)!}=
\frac{(n-k-l)!\,(n-j-l)!}{(n-\alpha+i-l+\rho)!}
\label{eq:ways}
\end{equation}
ways. 

Next, we apply the homomorphism $\iota$, see \eqref{eq:iota},
to the equation \eqref{eq:coefficients}.
By \eqref{eq:ways},
we get
$$
\iota\bigl(\Xi(\sigma)\bigr)\, \iota\bigl(\Xi(\tau)\bigr)
\cdot\Bigl( \frac{1}{n!}\cdot
 \frac{(n-k-l)!\,(n-j-l)!}{(n-\alpha+i-l+\rho)!}\Bigr)
=
\frc^\kappa_{\sigma\tau}(n)\, \iota\bigl(\Xi(\kappa)\bigr).
$$

By Lemma \ref{l:number-bijections}, for $\phi\in \PB(\alpha)$
$$
\iota\bigl(\Xi(\phi)\bigr)=\frac{n!}{(n-\alpha+\rk \phi)!}.
$$
We come to the equation
\begin{multline*}
\frac{n!}{(n-\alpha+i+j)!}\cdot\frac{n!}{(n-\alpha+i+k)!}
\cdot \frac{(n-k-l)!\,(n-j-l)!}{n!\,(n-\alpha+i-l+\rho)!}
= \\=\frc^\kappa_{\sigma\tau}(n)\cdot \frac{n!}{(n-\alpha+i+\rho)!}.
\end{multline*}
Keeping in mind \eqref{eq:ijkl} we come to
$$
\frc^\kappa_{\sigma\tau}(n)=
\frac{(n-\alpha+i+\rho)!}{(n-\alpha+i-l+\rho)!},
$$
i.e., to the desired formula \eqref{eq:constants}.

\sm

{\bf \punct 
Structure of algebras $\boldsymbol{\bbO[\alpha;n]}$.\label{ss:structure}}
%Dimensions of spaces of $\boldsymbol{S_n}$-fixed vectors. Proof of Proposition
%\ref{pr:isomorphic}.\label{ss:young}}
Recall that irreducible representations $\rho(N;\lambda)$
 of a symmetric group
$S_N$ are canonically enumerated by Young diagrams $\lambda$ of size
$|\lambda|=N$, see, e.g., \cite{Jam}. 
The dimension of $\rho(N;\lambda)$ is the number of Young tableaux
of the shape%
\footnote{There is a convenient `hook formula' for the dimensions but we do not need it.} $\lambda$. Recall that a  {\it Young   tableau}
of a shape $\lambda$ is a filing of boxes of $\lambda$
by $1$, \dots, $N$
such that numbers decrease in the right direction
and in the downward direction.
 
 The trivial (one-dimensional)
representation of $S_N$ corresponds to the row of length $N$.

\begin{proposition}
\label{pr:structure}
For all $n\ge\alpha$ algebras $\bbO[\alpha;n]$ are 
(noncanonically) isomorphic. Moreover, they are isomorphic
to the direct sum
\begin{equation}
\bigoplus_{t,\lambda'}\End(\C^{\ell(t,\lambda')})
\label{eq:direct-sum}
\end{equation}
 of matrix algebras of orders
\begin{equation}
%\Bigl\{
\ell(t,\lambda')= \frac{\alpha!}{t!\,(\alpha-t)!}\,
 \dim\rho(\alpha-t;\lambda')
%\Bigr
%\}, 
\label{eq:multi}
\end{equation}
where $t$ ranges in the set $\{0,1,\dots,\alpha\}$
and $\lambda'$ ranges in the set of Young diagrams of size
 $\alpha-t$.
\end{proposition}

{\bf\punct Proof of Proposition \ref{pr:structure}.}
Consider a representation $\rho(N;\lambda)$ of $S_N$
and its  restriction   to a subgroup $S_{N-k}$.
A representation $\rho(N-k;\mu)$ is contained in
the spectrum of the restriction iff the Young diagram $\mu$ is 
 contained in $\lambda$. 
The multiplicity of $\rho(N-k;\mu)$ in the restriction
 is the number of skew Young tableaux
of shape $\lambda\setminus\mu$. Recall that
 a {\it skew Young tableau} is a filing
of boxes of $\lambda\setminus \mu$ by $1$, \dots, $k$, such that 
 numbers decrease in the right direction
and in the downward direction.  

\sm

So, a representation $\rho(\alpha+n;\lambda)$ contains an $S_n$-fixed
vector if and only if the length of the  first row of $\lambda$
is $\ge n$. In this case $\mu$ consists of the first $n$ boxes 
of the first strip.
 If $n\ge\alpha$, then $\lambda\setminus \mu$
consists of two disjoint parts,  the rest $\kappa$ of the first
row and the union $\lambda'$ of the remaning rows.
For a filling of $\lambda\setminus \mu$, we must choose a collection of numbers
$1$, \dots, $\alpha$ for boxes of
 $\kappa$, and fill 
 $\lambda'$ by the remaining numbers. The latter question
is equivalent to a filling of Young tableau of the shape $\lambda'$.
We get the following statement.

\begin{proposition}
The collection of dimensions of nontrivial spaces of 
$S_n$-fixed vectors in irreducible representations of
$S_{n+\alpha}$ is \eqref{eq:multi}.
%\begin{equation}
%\Bigl\{ \frac{\alpha!}{t!\,(\alpha-t)!}\,
% \dim\rho(\alpha-t;\lambda')
%\Bigr\}, 
%\label{eq:multi}
%\end{equation}
%where $t$ ranges in the set $\{0,1,\dots,\alpha\}$
%and $\lambda'$ ranges in the set of Young diagrams of size
% $\alpha-t$.
\end{proposition}

This collection does not depend on $n$. Keeping in mind
\eqref{eq:oplus}, we get Proposition \ref{pr:structure}.

\section{The interpolation and the limit at $\infty$}

\COUNTERS

In this section, we construct a family of algebras
$\bbO[\alpha;\nu]$, where $\nu$ ranges in $\C$,
interpolating algebras $\bbO[\alpha;n]$ (Theorem \ref{th:interpolation},
we simply replace integer $n$ to complex $\nu$ in formulas
for structure constants). Next, we show that our formulas
imply an old Olshanski's
result \cite{Olsh} that the limit of such algebras
as $\nu\to\infty$ is the semigroup algebra of $\PB(\alpha)$
(Theorem \ref{th:limit}). Subsections \ref{ss:structure1},
\ref{ss:filtration}--\ref{ss:iota}
contain some simple observations on algebras $\bbO[\alpha;\nu]$.

\sm

{\bf \punct Varieties of associative algebras.} We need some preliminaries.
Consider a space $\C^M$ with a fixed basis $e_j$.
Consider a structure of an associative algebra on $\C^M$,
$$
e_p e_q= \sum c^s_{pq} e_s.
$$
Such product is associative if the {\it structure constants} 
$c^s_{pq}$ satisfy the following collection of $M^4$
 quadratic equations:
\begin{equation}
\sum_{i} c_{pq}^i c_{ir}^t=\sum_j c_{pj}^tc_{qr}^j.
\label{eq:ass}
\end{equation}
So the set $\Assoc(M)$  of all structures  
of associative algebras on $\C^M$
is an algebraic variety (\cite{Ger}). It has many irreducible components
of different dimensions  (see, e.g., \cite{Maz}, \cite{Ner-Lie}).
 The group $\GL(M,\C)$ acts on 
$\Assoc(M)$ by changes of bases. The set of its orbits is the set
of $M$-dimensional associative algebras defined up to isomorphisms.

 Consider a direct sum
of several matrix algebras, say, 
$$\cE(\{l_j\}):=\bigoplus_{j=1}^\tau \End(\C^{l_j}),$$
 such that
$\sum_{j=1}^\tau l_j^2=M$. 
 Consider the
subset $\cO(\{l_j\})$ in  $\mathrm{Assoc}(M)$ consisting of algebras
isomorphic to  $\cE(\{l_j\})$,
i.e., the $\GL(M,\C)$-orbit of $\cE(\{l_j\})$.
Clearly, 
\begin{equation}
\dim \cO(\{l_j\})=M^2-M+\tau.
\label{eq:MM}
\end{equation}
 Indeed, the stabilizer is the group of automorphsism of 
 $\cE(\{l_j\})$. Its connected component%
 \footnote{The group of automorphisms also contains permutations
 of factors if some $l_j$ are equal.}
 is the product of the groups of automorphisms of the summands,
 i.e., $\prod_j (\GL(l_j,\C)/\C^\times)$, where $\C^\times$ is
 the multiplicative group of $\C$ (embedded to $\GL(l_j,\C)$
 as scalar matrices). So, the dimension is $\sum (l_j^2-1)=M-\tau$.

 It  is well-known, see, e.g., \cite{Fla}%
 \footnote{We say, that a {\it deformation} of an algebra
 $\cA$ is a holomorphic curve $\cA(\epsilon)\subset \Assoc(M)$
 defined for small $|\epsilon|$
and satisfying $\cA(0)=\cA$. It is easy to verify, that if there is a sequence
 $\epsilon_k\to 0$ such that $\cA(\epsilon_k)$ have a nil-radical of dimension $l$, then the nil-radical of $\cA$ has dimension $\ge l$. So,
 deformations of semisimple algebras are semisimple. If for some 
 $\epsilon_k\to 0$ dimensions of the centers of $\cA(\epsilon_k)$
 are $p$, then the dimension of the center of $\cA$ is $\ge p$.
 If $\cA(0)\in \cO(\{l_j\})$ and the curve does not contained in this family, it must be contained in a family $\cO(\{l'_i\})$
 of larger dimension. By \eqref{eq:MM}, this requires a growth of a
 dimension of center, and we come to a contradiction.}, that each  closure
$\ov{\cO(\{l_j\})}\subset \Assoc(M)$ of 
$\cO(\{l_j\})$ is an irreducible component of the variety $\Assoc(M)$.
 
\sm

{\bf \punct The interpolation.}

\begin{theorem}
\label{th:interpolation}
Let $\nu\in\C$. For each triple $\sigma$, $\tau$, $\kappa\in\PB(\alpha)$
we  define constants
$\frc^\kappa_{\sigma,\tau}(\nu)$ in the following way:

\sm

--- if $\kappa$ is an extension of $\sigma\tau$, then
\begin{multline}
\frc^\kappa_{\sigma,\tau}(\nu)=\\=
(\nu-\alpha+\rk \kappa)(\nu-\alpha+\rk \kappa-1)\dots 
(\nu-2\alpha+\rk \sigma+\rk \tau-\rk \sigma\tau+\rk\kappa +1)
=\\=(-1)^l \bigl(-\nu+\alpha-\rk \kappa\bigr)_l,
\end{multline}
where 
$$l:=\alpha+\rk\sigma\tau-\rk\sigma-\rk\tau$$
and 
$$(p)_k=p(p+1)\dots (p+k-1)$$
 is the Pochhammer symbol.

\sm

--- 
$\frc^\kappa_{\sigma,\tau}(\nu)=0$ otherwise.

\sm

Consider a linear space $\C^{\# \PB(\alpha)}$,
 whose basis consists of symbols 
$\Xi(\sigma)$, where $\sigma$ ranges in $\PB(\alpha)$.
Then for any $\nu\in\C$ the formula
$$
\Xi(\sigma)\, \Xi(\tau)=\sum_\kappa \frc^\kappa_{\sigma,\tau}(\nu)
\,\Xi(\kappa) 
$$
determines a structure of an associative algebra on 
$\C^{\# \PB(\alpha)}$.
\end{theorem}

{\sc Proof.} We must show that for each $\nu$ the polynomials
$\frc^\kappa_{\sigma,\tau}(\nu)$ satisfy the system of quadratic equations
\eqref{eq:ass}. By Theorem \ref{th:coefficients}, the equality
is valid
for integer $\nu\ge\alpha$.
%So, this is a collection of polynomial conditions for $\nu$,
%which are valid for all integer $\nu\ge\alpha$. 
Therefore, they are valid for   all $\nu\in\C$.
\hfill $\square$

\sm

We denote such algebras by $\bbO[\alpha;\nu]$.

\sm

{\bf\punct Structure of the algebras $\boldsymbol{\bbO[\alpha;\nu]}$.%
\label{ss:structure1}}

\begin{proposition}
{\rm a)}
For all $\nu\in\C$ except a finite number of values
algebras $\bbO[\alpha;\nu]$ are (noncanonically).

\sm

{\rm b)} For the remaining exceptional values $\nu_p$ the algebras
$\bbO[\alpha;\nu_p]$ 
are not semisimple.
\end{proposition}

{\sc Remark.} Therefore, generic algebras $\bbO[\alpha;\nu]$ are isomorphic to algebras \eqref{eq:direct-sum}--\eqref{eq:multi}.
\hfill $\boxtimes$

\sm

{\sc Proof.} a) We have a reducible  algebraic variety
$\Assoc(M)$, where $M:=\#\PB(\alpha)$. We have the algebraic curve in
 $\Assoc(M)$
defined by $\nu\mapsto \bbO[\alpha;\nu]$. For $n\ge\alpha$
the point $\bbO[\alpha;n]$ is contained in a certain 
$\GL(M,\C)$-orbit
 $\cO(\{l_j\})$. Hence this curve is contained in the component 
  $\ov{\cO(\{l_j\})}$. The set 
  $\ov{\cO(\{l_j\})}\setminus \cO(\{l_j\})$
  has codimension $\ge 1$ in $\ov{\cO(\{l_j\})}$, and therefore only finite number of points $\nu_p$ of the curve are not contained
  in $\ov{\cO(\{l_j\})}$.
  
  \sm
  
  b) If for an exceptional point $\nu_p$ the algebra
   $\bbO[\alpha;\nu_p]$
  is semisimple,  then it is a point of some orbit  $\cO(\{m_p\})$.
  So the  whole curve $\nu\mapsto \bbO[\alpha;\nu]$ except
  a finite number of points is contained in the same
   $\GL(M,\C)$-orbit, and
   $\bbO[\alpha;\nu_p]$ is not an exception.  
  \hfill $\square$

\sm

{\bf \punct The limit as $\boldsymbol{\nu\to\infty}$.}
Consider another normalization of the natural basis in
$\bbO[\alpha;\nu]$,
$$
\xi_\nu(\sigma):= \frac1{a(\nu;\sigma)}\,
\Xi(\sigma),
$$
 where
 \begin{equation}
  a(\nu;\sigma):=\nu(\nu-1)\dots(\nu-\alpha+\rk\theta+1).
  \label{eq:anu}
\end{equation}
So, for integer $\nu=n\ge\alpha$ we have
$$
\iota(\xi(\sigma))=1,
$$
where $\iota(\cdot)$ is defined by \eqref{eq:iota}.

Consider the structure constants of our algebra in the new basis,
$$
\xi(\sigma)\, \xi(\tau)=\sum_\kappa \frd^\kappa_{\sigma,\tau}(\nu)
\,
\xi(\kappa),
$$
where
$$
 \frd^\kappa_{\sigma,\tau}= \frc^\kappa_{\sigma,\tau}
 \cdot \frac {a(\nu;\kappa)}{a(\nu;\sigma)\,a(\nu;\tau)}.
$$
We have the following asymptotics
as $|\nu|\to\infty$:
\begin{align*}
a(\nu;\phi)&\sim \nu^{\alpha-\rk\phi};\\
\qquad \frc^\kappa_{\sigma,\tau}&\sim
\nu^{\alpha-\rk\sigma-\rk\tau+\rk \sigma\tau}
\quad\text{if $\kappa$ is an extension of $\sigma\tau$.}
\end{align*}
Therefore,
$$
 \frd^\kappa_{\sigma,\tau}
 \sim\frac{\nu^{\alpha-\rk\sigma-\rk\tau+\rk \sigma\tau}
\nu^{\alpha-\rk \kappa} }{\nu^{\alpha-\rk\sigma}\nu^{\alpha-\rk\tau}}
=\nu^{\rk\sigma\tau-\rk\kappa}.
$$
Since $\kappa\sqsupset\sigma\tau$, we get
$$
\lim_{|\nu|\to\infty}\frd^\kappa_{\sigma,\tau}=
\begin{cases}
1,\qquad \text{if $\kappa=\sigma\tau$};
\\
0,\qquad \text{otherwise.}
\end{cases}
$$

Thus we get the following theorem (Olshanski \cite{Olsh1}, 
\cite{Olsh}).

\begin{theorem}
\label{th:limit}
The natural limit  of $\bbO[\alpha;\nu]$ 
as $\nu\to \infty$ is isomorphic to 
the semigroup algebra%
\footnote{Let $\Gamma$ be a finite semigroup. The {\it semigroup
algebra} $\C[\Gamma]$ has a basis consisting of symbols $\delta_g$
and the product is defined by $\delta_{g_1}\cdot\delta_{g_2}=\delta_{g_1g_2}$.} of $\PB(\alpha)$.
\end{theorem}

\begin{proposition}
\label{pr:lim-isomorphic}
The semigroup algebra $\C[\PB(\alpha)]$ of $\PB(\alpha)$
is (noncanonically) isomorphic to algebras $\bbO[\alpha;n]$
for $n\ge\alpha$. 
\end{proposition}

Recall a general statement. The semigroup algebra $\C[\Gamma]$
 of a finite 
inverse semigroup% $\Gamma$
\footnote{An {\it inverse semigroup} is a subsemigroup of 
a semigroup of partial bijections closed with respect to the pseudoinversion.}
is semisimple. It is the sum of matrix algebras 
$$
\bigoplus_{\rho\in \wh\Gamma} \End(V_\rho),
$$
where the summation is taken over all irreducible representations
of $\Gamma$; $V_\rho$ are  spaces
of representations (see Oganesyan \cite{Oga} or Clifford, Preston
 \cite{CP}, Theorem 5.26).

\sm

{\sc Proof of Proposition \ref{pr:lim-isomorphic}.} 
We apply this statement to the semigroup algebra $\C[\PB(\alpha)]$.
Irreducible representations of $\PB(\alpha)$
 can be easily described, see, e.g., \cite{Ner-book},
Subsect.~VIII.2.2. They are enumerated by
pairs $(t, \mu)$, where  $t=0$, 1, \dots, $\alpha$,
and $\mu$ is a Young diagram of size $\alpha-t$.
The explicit description shows that such representation
is realized in a space of functions on the set
of $(\alpha-t)$-element subsets in $\Omega_\alpha$
taking values in the space of representation $\rho(\alpha-t;\mu)$
of $S_{\alpha-t}$ (we keep the notation of the Subsect. \ref{ss:structure}).
So, the dimension of such a representation
is
$$
\frac{\alpha!}{t!\,(\alpha-t)!}\, \dim\rho(\alpha-t;\mu),
$$
and we get the same collection of dimensions \eqref{eq:multi}.
\hfill $\square$

\sm

{\bf \punct The filtration in 
$\boldsymbol{\bbO[\alpha;\nu]}$.%
\label{ss:filtration}} 
Let $s=0$, $1$, \dots, $\alpha$. 
Denote by $\bbO_s[\alpha;\nu]$ 
the linear subspace in $\bbO_s[\alpha;\nu]$
generated by elements $\Xi(\sigma)$ such that
$\rk\sigma\ge \alpha-s$. So, 
$$
\bbO_0[\alpha;\nu]\subset \bbO_1[\alpha;\nu]\subset
\dots \subset \bbO_\alpha[\alpha;\nu],
$$
and
$$\bbO_0[\alpha;\nu]=\C[S_\alpha],\qquad
 \bbO_\alpha[\alpha;\nu]=\bbO[\alpha;\nu].$$

\begin{proposition}
\label{pr:filtration}
{\rm a)}
Let $s+t\le\alpha$.
Let $\Xi(\sigma)\in \bbO_s[\alpha;\nu]$, 
$\Xi(\tau)\in\bbO_t[\alpha;\nu]$. Then
 $$\Xi(\sigma)\,\Xi(\tau)-\Xi(\sigma\tau)\in \bbO_{s+t-1}
 [\alpha;\nu].$$
 
{\rm b)} We have
$$\bbO_s[\alpha;\nu]\cdot \bbO_t[\alpha;\nu]=
\begin{cases}\bbO_{s+t}[\alpha;\nu]&\quad \text{if $s+t<\alpha;$}
\\
\bbO_\alpha[\alpha;\nu]&\quad\text{otherwise}
.
\end{cases}
$$

\sm

{\rm c)} The subspace $\bbO_1[\alpha;n]$ generates the whole algebra
$\bbO[\alpha;\nu]$.
  
\end{proposition}

{\sc Proof.} a) The statement  immediately follows from
the formula for structure constants.  

\sm

b) The inclusion $\subset$ follows from a). Let us prove $\supset$.
Suppose that the statement is valid if $s+t\le j-1$.
 Let us show that this is so 
 if $s+t=j$. For each $\mu\in\PB(\alpha)$
  of rank $n-j$ we can find 
 $\sigma$ of rank $n-s$ and $\tau$ of rank $n-t$ such that
 $\sigma\tau=\mu$. Then
 $$
 h:=\Xi(\sigma)\,\Xi(\tau)-\Xi(\mu)\in \bbO_{j-1}(\alpha;\nu).
 $$
Since $\Xi(\sigma)\,\Xi(\tau)\in \bbO_{j}(\alpha;\nu)$ and $h\in \bbO_{j-1}(\alpha;\nu)$, we get $\Xi(\mu)\in \bbO_{j}(\alpha;\nu)$.

\sm

c) This follows from b).
\hfill $\square$

\sm

{\bf \punct The trace.}
We define the {\it trace} $\Tr(\cdot):\bbO[\alpha;\nu]\to\C$
  as a linear functional
%$\Tr:\bbO[\alpha;\nu]\to\C$
 defined on basis elements \eqref{eq:basis}
by
\begin{align*}
\Tr\bigl(\Xi(1))&=1\\
\Tr(\Xi(\sigma))&=0\quad \text{for $\sigma\ne 1$.}
\end{align*}

\begin{proposition}
\label{pr:trace}
\begin{equation}
\Tr(H_1 H_2)=\Tr(H_2 H_1).
\end{equation}
\end{proposition}

{\sc Proof.} For $\nu=n\ge\alpha$ we have
$$
\Tr H=n! \tr_ H,
$$
where $\tr(\cdot)$ is the trace  \eqref{eq:trace}
on $\C[S_n\backslash S_{\alpha+n}/S_n]$. Hence 
the collection of equalities 
$$
\Tr\bigl(\Xi(\sigma)\, \Xi(\tau) \bigr)
=\Tr\bigl( \Xi(\tau)\, \Xi(\sigma) \bigr)
$$
 is valid for such $\nu$. Since both sides are polynomials in $\nu$,
 this is valid for all complex  $\nu$.
\hfill $\square$
 
 \sm

{\bf \punct The involution and the inner product.} Now let $\nu\in \R$.
Consider the antilinear involution of
$\bbO[\alpha;\nu]$ 
defined by
$$
\Xi[\sigma]^*=\Xi[\sigma^\square].
$$
For $\nu\in\R$ 
we define the Hermitian form on $\bbO[\alpha;\nu]$
by 
$$
\la H_1,H_2\ra:=\Tr H_1 H_2^*.
$$

\begin{proposition}
\label{pr:positivity}
For sufficiently large $\nu\in \R$, this form is positive definite.
\end{proposition}

{\sc Proof.}  The positivity of our form is equivalent to a
 finite collection of polynomial inequalities on $\nu$. For instance, we can require the
 positivity of all principal minors of the matrix
 $$
 \bigl\{\Tr\bigl(\Xi(\sigma)\,\Xi(\tau^\square)\bigr)\bigr\}_{\sigma,\tau\in \PB(\alpha)}.
 $$
  These inequalities take place for all integer $\nu \ge\alpha$,
 therefore they are valid for sufficiently large $\nu$.
 \hfill $\square$ 

\sm

{\bf \punct The homomorphism $\boldsymbol{\iota}$.%
\label{ss:iota}}
Consider the linear functional $\iota$ on 
$\bbO[\alpha;\nu]$
defined on basis elements by
$$
\iota(\Xi(\sigma))=a(\nu;\sigma),
$$
where $a(\nu;\sigma)$ is given by \eqref{eq:anu}.

\begin{proposition}
$\iota$  is a homomorphism $\bbO[\alpha;\nu]\to\C.$
\end{proposition}

{\sc Proof.}  
Indeed, for $\nu=n\ge\alpha$ the  functional $\iota$ coincides 
with the homomorphism $\iota:\bbO[\alpha;n]\to\C$ (see \eqref{eq:iota}).
So, for any pair $\sigma$, $\nu\in \PB(\alpha)$
 the equality
$$\iota\bigl(\Xi(\sigma)\bigr)\,\iota\bigl(\Xi(\tau)\bigr)=
\iota\bigl(\Xi(\sigma)\, \Xi(\tau)\bigr)$$
is valid for such $\nu$. Since this condition is 
 an equality for polynomials in  $\nu$,
  it is valid for all $\nu\in\C$.
 \hfill $\square$

\section{Generators and relations}

\COUNTERS

Now we  describe algebras $\bbO[\alpha;\nu]$ in terms of generators and relations (Theorem \ref{th:homomorphism}).

\sm

{\bf \punct Generators and relations.} 
%Denote by $(ij)\in S_\alpha$ the transposition of $i$ and $j$.
%
 Let $\alpha\ge 0$ be integer. 
Let $\nu\in\C$.
 Consider the associative algebra $\wt\bbO[\alpha,\nu]$
with generators 
$$
\text{$A(g)$, where $g$ ranges in $S_\alpha$,
and $\Theta_1$, \dots, $\Theta_\alpha$}
$$
and the relations
\begin{align}
A(g_1)A(g_2)&=A(g_1 g_2);
\label{eq:rel1}
\\
A(g)\Theta_i A(g^{-1})&=\Theta_{g(i)};
\label{eq:rel2}
\\
\Theta_j^2&=(\nu-1)\Theta_j+ \nu A(1);
\label{eq:rel3}
\\
\Theta_j \Theta_i-\Theta_i \Theta_j&=A((ij))\,(\Theta_j-\Theta_i);
\label{eq:rel4}
\\
\bigl(A((ij))-A(1)\bigr)\, \Theta_j (\Theta_i+A(1))&=0.
\label{eq:rel5}
\end{align}

%\sm

%{\sc Remark.} So $A(1)$ is a unit of the algebra, 
%and all elements $A(g)$ are invertible, $A(g)^{-1}=A(g^{-1})$. 
%\hfill $\boxtimes$

 For any $i\le \alpha$ we  denote the following elements 
 of $\bbO[n;\nu]$:
 $$
 \theta_i:=
 \Xi\bigr(T\{1,\dots,i-1,i+1,\dots,\alpha\}\bigr)\in \bbO[\alpha;\nu]
 =: \Xi\bigr(T\{1,\dots,\wh i,\dots,\alpha\}\bigr)
 .
 $$

\begin{theorem}
\label{th:homomorphism}
%{\rm a)} $\dim  \bbO[\alpha,\nu]=\#\Pi_\alpha$.
%
%\sm 
%
% Let $n\ge 0$ be integer. 
 Then the map
 $$
 \Theta_i\mapsto \theta_i, \qquad A(g)\mapsto \Xi(g) 
 $$
 defined on the generators 
 determines an isomorphism $\wt\bbO[\alpha;\nu]\to \bbO[\alpha;\nu]$.
 
% \sm 
 
% {\rm c)} If $n\ge \alpha$, then this homomorphism is an isomorphism.
% If $0\le n<\alpha$, then it is surjective.
 
%\sm  
 
% {\rm d)} For all $\nu\in \C$ except finite collection
% of values, the algebras $\bbO[\alpha;\nu]$ are%
% \footnote{Noncanonically.}  isomorphic
% (and are direct sums of matrix algebras). 
 \end{theorem}
 
 The proof occupies the rest of this section.
 
 \sm
 
{\bf \punct Relations.}

\begin{lemma}
If we substitute $A(g)\mapsto\Xi(g)$ and $\theta_i\mapsto \Theta_i$
to the relations \eqref{eq:rel1}--\eqref{eq:rel5},
then we get correct correct identities.
\end{lemma} 

{\sc Proof.}  For integer  $\nu=n\ge\alpha$
the relations \eqref{eq:rel1}--\eqref{eq:rel2}
are obvious. The relation \eqref{eq:rel3} is a special case
of Theorem \eqref{th:coefficients}.

Denote
$$
\theta_{ij}=\theta_{ji}:=
\Xi\bigl(T\{1,\dots,\wh i,\dots,\wh j,\dots,\alpha\}\bigr).
$$
Multiplying $\theta_i\cdot\theta_j$ with Theorem \ref{th:coefficients},
we get
$$
\theta_i \theta_j=\theta_{ij}+\Xi((ij))\theta_j.
$$
Subtracting $\theta_i\theta_j-\theta_j\theta_i$
we come to the relation  \eqref{eq:rel4}.

It remains to verify the relation \eqref{eq:rel5}.
Clearly, 
$$
\Xi((ij))\theta_{ij}=\theta_{ij}.
$$
Therefore, 
\begin{multline*}
0=\bigl(\Xi((ij))-1\bigr)\theta_{ij}=\bigl(\Xi((ij))-1\bigr)
\bigl(\theta_j\theta_i-\Xi((ij))\theta_j\bigr)
=\\=
\bigl(\Xi((ij))-1\bigr)\theta_j\theta_i-\bigl(\Xi((ij))-1\bigr)\Xi((ij))\theta_j
=\\=\bigl(\Xi((ij))-1\bigr)\theta_j\theta_i+\bigl(\Xi((ij))-1\bigr)\theta_j
=\bigl(\Xi((ij))-1\bigr)\theta_j(\theta_i+1).
\end{multline*}

Since the relations are valid for integer $\nu\ge \alpha$,
they are valid for all $\nu$. \hfill $\square$

\sm 

Thus, the homomorphism 
$\wt \bbO[\alpha;\nu]\to\bbO[\alpha;\nu]$ is well defined.

\sm

{\bf\punct The surjectivity of the homomorphism 
$\wt \bbO[\alpha;\nu]\to\bbO[\alpha;\nu]$.}
%We refer to Proposition \ref{pr:filtration}.
Any partial bijection $\sigma$ of rank $\alpha-1$
 can be represented
as a product of the form $g\cdot T\{1,\dots, \wh i,\dots,\alpha\}$,
see \eqref{eq:gT}.
Therefore, $\Xi[\sigma]$ is contained in the image of the homomorphism.
Applying Proposition \ref{pr:filtration}.c, we get the desired surjectivity.

\begin{corollary}
\label{cor:ge}
$\dim \wt \bbO[\alpha;\nu]\ge \dim  \bbO[\alpha;\nu]=\#\PB(\alpha)$.
\end{corollary}
 
{\bf \punct A basis in $\boldsymbol{\wt \bbO[\alpha;\nu]}$.}

\begin{theorem}
\label{th:basis}
The following elements form a basis in the algebra
$\wt \bbO[\alpha;\nu]$:
\begin{multline} \label{eq:basis}
A(g)\, \Theta_{i_1}\dots \Theta_{i_k},
\qquad\text{where $i_1<i_2<\dots<i_k$ and}\\  g(i_1)<g(i_2)<\dots<g(i_k).
\end{multline}
\end{theorem}

The statement is proved in the next subsection.

\sm

{\bf \punct The dimension of $\boldsymbol{\wt \bbO[\alpha;\nu]}$.}

\begin{lemma}
\label{l:span-of-basis}
The linear span of elements \eqref{eq:basis}
is the whole $\wt \bbO[\alpha;\nu]$.
\end{lemma}

{\sc Proof.} Denote by $\wt \bbO_t[\alpha;\nu]$
the subspace in $\wt \bbO[\alpha;\nu]$
spanned by  all products of the form
$$
A(g_1)\,\Theta_{i_1} A(g_2)\,\Theta_{i_2}
\dots A(g_p)\,\Theta_{i_p} A(g_{p+1}),
\qquad\text{where $p\le t$.}
$$

Clearly,
$$
\wt \bbO_s[\alpha;\nu]\cdot \wt \bbO_t[\alpha;\nu]
=
\wt \bbO_{s+t}[\alpha;\nu].
$$

By the relations \eqref{eq:rel1}-\eqref{eq:rel2}, we get:

\sm

--- {\it Any subspace $\wt \bbO_t[\alpha;\nu]$
is spanned by products  of the form
\begin{equation}
A(g)\Theta_{i_1}\Theta_{i_2}
\dots \Theta_{i_p},
\qquad\text{where $p\le t$.}
\label{eq:spanned}
\end{equation}
}

Next, for any permutation $\lambda$ of a set $\{1,\dots,t\}$
we have 
\begin{equation}
\Theta_{i_{\lambda(1)}}\dots \Theta_{i_{\lambda(t)}}
- \Theta_{i_1}\dots \Theta_{i_t}\in \wt \bbO_{t-1}[\alpha;\nu]
\label{eq:TTTT}
\end{equation}
(we apply \eqref{eq:rel4} many times).
%Indeed, we decompose $\lambda$ as a product of transpositions
%and many times apply \eqref{eq:rel2}. 

Therefore {\it 
  products  \eqref{eq:spanned}
with $i_1\le i_2\le \dots \le i_p$ span $\bbO_t[\alpha;\nu]$}.
We verify this claim by induction. For $t=1$ this is correct. 
We assume that this is correct until $t-1$, and 
 \eqref{eq:TTTT} proves the induction step.

 By \eqref{eq:rel2},
 $\Theta_j^2\in  \wt \bbO_1[\alpha;\nu]$.
 The same induction arguments imply the following statement:
 
\sm 
 
 --- {\it Any subspace $\wt \bbO_t[\alpha;\nu]$
is spanned by products  of the form
\begin{equation}
A(g)\Theta_{i_1}\Theta_{i_2}
\dots \Theta_{i_p},
\qquad\text{where $p\le t$ and $i_1<i_2<\dots<i_p$.}
\label{eq:spanned1}
\end{equation}
}

\sm

In particular, $\wt\bbO_\alpha[\alpha;\nu]=\wt\bbO[\alpha;\nu]$.

\begin{lemma}
\label{l:independence}
Let  $h\in S_\alpha$ be supported by a set $i_1<\dots<i_p$.
Then 
$$
A(h)\Theta_{i_1}\Theta_{i_2}
\dots \Theta_{i_p}-\Theta_{i_1}\Theta_{i_2}
\dots \Theta_{i_p}\in \wt \bbO_{p-1}[\alpha;\nu].
$$
\end{lemma}

{\sc Proof of Lemma \ref{l:independence}.}
Representing the relation \eqref{eq:rel5}
in the form
\begin{equation*}
A((ij))\, \Theta_j\Theta_i-\Theta_j\Theta_i=\Theta_j- A((ij))\Theta_j,
\end{equation*}
we see 
\begin{equation}
\label{eq:erase} 
 A((ij))\, \Theta_j\Theta_i-\Theta_j\Theta_i\in \wt\bbO_1[\alpha;\nu].
\end{equation}

Decompose our $h$ in a product of transpositions of 
$i_1$, \dots, $i_p$, say
$$h=(i_{u_1}i_{v_1})\dots (i_{u_m}i_{v_m}).$$
So,
$$
A(h)\, \Theta_{i_1}\dots \Theta_{i_p}= 
A\bigl((i_{u_1}i_{v_1})\bigr)\dots A\bigl((i_{u_m}i_{v_m})\bigr)\,
\Theta_{i_1}\dots \Theta_{i_p}.
$$
Transposing factors $\Theta$ we move $\Theta_{i_{u_m}}$ to the first position, each such transposition changes the expression by a summand
$\in \wt\bbO_{p-1}[\alpha;\nu]$, see \eqref{eq:TTTT}.
Then we  move
$\Theta_{i_{v_m}}$ to the second position.
Using \eqref{eq:erase}, 
%$$
%\deg_\Theta\Bigl(A\bigl((i_{p_m}i_{q_m})\bigr) \Theta_{i_{p_m}}\Theta_{i_{q_m}}
%-\Theta_{i_{p_m}}\Theta_{i_{q_m}}\Bigr)
%<2.
%$$
 we erase $A\bigl((i_{u_m}i_{v_m})\bigr)$; afterwards we return $\Theta_{i_{u_m}}$ and 
$\Theta_{i_{v_m}}$ to the initial positions. We get
\begin{multline*}
\Bigl(A\bigl((i_{u_1}i_{v_1})\bigr)\dots A\bigl((i_{u_{m-1}}i_{v_{m-1}})\,A\bigl((i_{u_m}i_{v_m})\bigr)-
\\-
A\bigl((i_{u_1}i_{v_1})\bigr)\dots A\bigl((i_{u_{m-1}}i_{v_{m-1}})\bigr)
\Bigr)\,
\Theta_{i_1}\dots \Theta_{i_p}\in \wt\bbO_{p-1}[\alpha;\nu].
\end{multline*}
Repeating this transformation, we come to the desired statement.
\hfill $\square$

\sm

{\sc The end of the proof of Lemma \ref{l:span-of-basis}.}
We prove the lemma by induction in the form:

\sm

--- {\it for each $s$, elements of the form
\begin{multline}
A(g)\,\Theta_{i_1}\Theta_{i_2}
\dots \Theta_{i_p},
\qquad\text{where $p\le s$ and 
$i_1<i_2<\dots<i_p$},\\ g(i_1)<\dots<g(i_p),
\label{eq:spanned2}
\end{multline}
generate $\bbO_s[\alpha;\nu]$.
}

\sm

Assume that this is valid for $s<t$. We take an element
\eqref{eq:spanned1} with $p=t$ and choose
$h\in S_\alpha$ supported by the set $i_1<\dots<i_t$ such that
$gh(i_1)<\dots<gh(i_t)$. Such $h$ is unique,
$$
h=h(i_1,\dots, i_t;g).
$$
By Lemma \ref{l:independence}, we have
$$
A(gh)\,\Theta_{i_1}\Theta_{i_2}
\dots \Theta_{i_t}-A(g)\,\Theta_{i_1}\Theta_{i_2}
\dots \Theta_{i_t}\in \bbO_{t-1}(\alpha;\nu).
$$ 
By the induction hypothesis, the difference is spanned
by elements of the desired type, and we can omit all
$A(g)\,\Theta_{i_1}
\dots \Theta_{i_t}$ such that $h(i_1,\dots, i_t;g)\ne 1$.

\begin{corollary}
\label{cor:le}
$\dim \wt\bbO[\alpha;\nu]\le \#\PB(\alpha)$.
\end{corollary} 

{\sc Proof.} For an expression \eqref{eq:spanned2},
we assign a partial bijection $\sigma\in\PB(\alpha)$
by the following rule:

\sm

--- $\dom\sigma$ is the complement to the set
$\{i_1, i_2,\dots, i_p\}$;

\sm

--- for $j\in \dom\sigma$ we assume $\sigma(j)=g(j)$.

\sm

Clearly, we get a bijection between the set of 
elements \eqref{eq:spanned2} (whose linear span is the whole algebra)
and $\PB(\alpha)$.
\hfill $\square$

\sm

{\sc Proofs of Theorems \ref{th:homomorphism} and \ref{th:basis}.}
Comparing Corollaries \ref{cor:ge} and \ref{cor:le}, we get
$$
\dim\wt\bbO[\alpha;n]=\#\PB(\alpha).
$$
Therefore, the elements \eqref{eq:basis} form a basis, and the
surjective homomorphism
$\wt\bbO[\alpha;\nu]\to \bbO[\alpha;\nu]$
is an isomorphism.
\hfill $\square$

%%%%%%%%%%%%%%%%%%%%%%%%%%%%%%%%%%%%%%%%%%%%
%%%%%%%%%%%%%%%%%%%%%%%%%%%%%%%%%%%%%%%%%%%%
%%%%%%%%%%%%%%%%%%%%%%%%%%%%%%%%%%%%%%%%%%%%
%%%%%%%%%%%%%%%%%%%%%%%%%%%%%%%%%%%%%%%%%%%%
%%%%%%%%%%%%%%%%%%%%%%%%%%%%%%%%%%%%%%%%%%%%

\noindent
\tt
University of Graz,
\\
\vphantom\hfill
Department of Mathematics and Scientific computing;
\\
High School of Modern Mathematics MIPT;
%\\
%Moscow State University, MechMath. Dept;
\\
University of Vienna, Faculty of Mathematics.
\\
e-mail:yurii.neretin(dog)univie.ac.at
\\
URL: https://www.mat.univie.ac.at/$\sim$neretin/
\\
\phantom{URL:}
https://imsc.uni-graz.at/neretin/index.html

\end{document}